\documentclass[a4paper,11pt]{article}
\usepackage[T2A]{fontenc} 

\usepackage[cp866]{inputenc} 

\usepackage[english]{babel} 
\usepackage{amssymb,amsmath,color} 

\newcommand{\eps}{\varepsilon}

\newcommand{\R}{{\mathbb{R}}}

\begin{document}

\newcounter{items}
\setcounter{items}{0}

\large \pagestyle{empty}
\begin{center}
 \quad UDC 517.925.42 \hfill MATHEMATICS \quad \ \  \\
 \bf ON MITROPOL'SKII YU.A.'S THEOREM ON PERIODIC SOLUTIONS OF
 SYSTEMS OF NONLINEAR DIFFERENTIAL EQUATIONS WITH
 NON-DIFFERENTIABLE RIGHT-HAND-SIDES\\
\end{center}
\centerline{ \bf  A. Buic\u{a}, J. Llibre, O. Yu. Makarenkov}

\

In the present paper we study the existence, uniqueness and
asymptotic stability of the $T$--periodic solutions for the system
\begin{equation}\label{ps}
\dot x =\eps g(t,x,\eps),
\end{equation}
where $\eps>0$ is a small parameter and the function $g\in
C^0(\mathbb{R} \times\mathbb{R}^k\times[0,1],\mathbb{R}^k)$ is
$T$--periodic in the first variable and locally  Lipschitz with
respect to the second one. As usual a key role will be played by
the averaging function
\begin{equation} \label{averfunction}
{g}_0(v)=\int\limits_0^T g(\tau,v,0)d\tau,
\end{equation}
and we shall look for those periodic solutions  that starts near
some $v_0\in {g}_0^{-1}(0)$.

In the case that $g$ is of class $C^1$, we remind the periodic
case of the second Bogolyubov's theorem (\cite{bog}, Ch.~1, \S~5,
Theorem~II) which represents a part of the averaging principle:
{\it $\det \, (g_0)'(v_0)\neq 0$ assures the existence and
uniqueness, for $\eps>0$ small, of a $T$--periodic solution of
system (\ref{ps}) in a neighborhood of $v_0,$ while the fact that
all the eigenvalues of the matrix $(g_0)'(v_0)$ have negative real
part, provides also its asymptotic stability.}

It was Mitropol'skii  who noticed that various applications
require the generalization of the second Bogolubov's theorem for
Lipschitz right hand parts. Assuming that $g$ is Lipschitz,
$g_0\in C^3(\R^k,\R^k)$ and that all the eigenvalues of the matrix
$(g_0)'(v_0)$ have negative real part Mitropol'skii developed an
analog of the second Bogolyubov's Theorem proving the existence
and uniqueness of a $T$--periodic solution of system (\ref{ps}) in
a neighborhood of $v_0.$ There was a great progress weakening the
assumptions of Mitropol'skii in his existence result (see
Samoylenko \cite{sam} and Mawhin \cite{maw1}) and some analogs of
his uniqueness result for equations with monotone nonlinearities
have been proposed (see papers by A.I. Perov, Yu.V. Trubnikov,
V.L. Hackevich \cite{per1}, \cite{per2}). Nevertheless, the
asymptotic stability conclusion of the second Bogolyubov's Theorem
remained to be not generalized for Mitropol'skii's settings
(namely, when $g$ is Lipschitz) for a long time. It has been done
recently by Buic\u{a}--Daniilidis in \cite{adr} for a class of
functions $v\mapsto g(t,v,0)$ differentiable at $v_0$ for almost
any $t\in[0,T],$ but it is assumed in \cite{adr} that the
eigenvectors of the matrix $({g}_0)'(v_0)$ are orthogonal and that
the function $g$ has a continuous Clarke differential, that is not
easy to check in applications.

In the next section of the paper assuming that $g$ is piecewise
differentiable in the second variable we show in Theorem~2 that
Mitropol'skii's conditions imply not only uniqueness, but also
asymptotic stability of a $T$--periodic solution of system
(\ref{ps}) in a neighborhood of $v_0.$ In other words we show that
Bogolyubov's theorem formulated above is valid when $g$ is  not
necessary $C^1$. Theorem~2 follows from our even more general
Theorem~1 whose hypotheses do not use any differentiability
neither of $g$ nor of $g_0$.

1.  Throughout the paper $\Omega\subset\R^k$ is some open set and
for any $\delta>0$ the $\delta$-neighborhood of $v\in\mathbb{R}^k$
is denoted by $B_\delta(v_0)=\left\{ v\in \R^k~:~\|v-v_0\|\leq
\delta \right\}$. We have the following main result.

{\bf Theorem 1.} {\it Let $g\in C^0(\mathbb{R}\times \Omega \times
[0,1],\mathbb{R}^k)$ and $v_0\in \Omega$. Assume the following
four conditions.
\begin{enumerate}
\setcounter{enumi}{\value{items}}
\item\label{lip} For some $L>0$ we have that
$\left\|g(t,v_1,\eps)-g(t,v_2,\eps)\right\|\le
L\left\|v_1-v_2\right\|$ for any $t\in[0,T],\ v_1,v_2\in \Omega,$
$\eps\in[0,1].$ \setcounter{items}{\value{enumi}}
\item\label{newunif} For any $\gamma>0$ there exists $\delta>0$ such that
\[
\begin{array}{l}
\hskip-1cm\left\|\int_0^T g(\tau,v_1+u(\tau),\eps)d\tau-\int_0^T
g(\tau,v_2+u(\tau),\eps)d\tau\right.\\
\left.-\int_0^T g(\tau,v_1,0)d\tau+\int_0^T
g(\tau,v_2,0)d\tau\right\|\le\gamma\|v_1-v_2\|
\end{array}
\]
for any $u\in C^0([0,T],\mathbb{R}^k),$ $\|u\|\le\delta,$
$v_1,v_2\in B_\delta(v_0)$ and $\eps\in[0,\delta].$
\setcounter{items}{\value{enumi}}
\item\label{def}Let $g_0$ be the averaging function given by
(\ref{averfunction}) and consider that ${g}_0(v_0)=0.$
\setcounter{items}{\value{enumi}}
\item\label{conde}
There exist $q\in[0,1),$ $\alpha,\delta_0>0$ and a norm
$\|\cdot\|_0$ on $\R^k$  such that $\left\|v_1+ \alpha
{g}_0(v_1)\right.$ $\left. -v_2- \alpha {g}_0(v_2)\right\|_0 \leq
q \|v_1-v_2\|_0$ for any $v_1,v_2\in  B_{\delta_0}(v_0).$
\setcounter{items}{\value{enumi}}
\end{enumerate}
Then there exists $\delta_1>0$ such that for every
$\eps\in(0,\delta_1]$ system (\ref{ps}) has exactly one
$T$--periodic solution $x_\eps$ with $x_\eps(0)\in
B_{\delta_1}(v_0).$  Moreover the solution $x_\eps$ is
asymptotically stable and $x_\eps(0)\to v_0$ as $\eps\to 0.$}

When solution $x(\cdot,v,\eps)$ of system (\ref{ps}) with initial
condition $x(0,v,\eps)=v$ is well defined on $[0,T]$ for any $v\in
B_{\delta_0}(v_0)$, the map $v\mapsto x(T,v,\eps)$ is well defined
and it is said to be the {\it Poincar\'e map} of system
(\ref{ps}). The proof of existence, uniqueness and stability of
the $T$--periodic solutions of system (\ref{ps}) in Theorem~1
reduces to the study of corresponding properties of the fixed
points of this map. More precisely, in order to prove Theorem~1 we
 represent $x(T,v,\eps)$  as
\[
x(T,v,\eps)=v+\eps g_\eps(v),\ \ {\rm where }\ \
g_\eps(v)=\int\limits_0^T g(\tau,x(\tau,v,\eps),\eps)d\tau.
\]
Then we show that the function $g_\eps$ satisfies a Lipschitz
condition uniformly in $\eps>0$ in the ball $B_{\delta_0}(v_0),$
where $\delta_0>0$ is a fixed small constant and moreover that
$g_\eps$ satisfies the following analog of the property (ii): {\it
for any  $\gamma>0$ there exists $\delta\in[0,\delta_0]$ such that
\[
\| g_\eps(v_1)-g_0(v_1)-g_\eps(v_2)+g_0(v_2)\| \le \gamma \|
v_1-v_2 \| \] for all  $v_1,v_2\in B_\delta(v_0)$ and
$\eps\in[0,\delta].$}

This allows to conclude that if $I+\eps g_0$ does not exceed
$1-\eps \widetilde{q}$ (where $\widetilde{q}$ is a constant) then
the Lipschitz constant of the function  $I+\eps g_\eps$ does not
exceed $1$ for $\eps>0$ sufficiently small. Thus the conclusion of
theorem~1 follows by applying the asymptotic stability of the
Poincar\'e map fixed point lemma (see e.g. \cite{kraop},
lemma~9.2). It turns out that the Lipschitz constant of the
function $I+\eps g_0$ does not exceed  $1-\eps \widetilde{q}$
indeed, where $\widetilde{q}=(1-q)/ \alpha,$ providing that the
assumptions of theorem~1 are satisfied.

In general it is not easy to check assumptions (\ref{newunif}) and
(\ref{conde}) in the applications of Theorem~1. Thus we give also
the following theorem based on Theorem~1 which assumes certain
type of piecewise differentiability instead of (\ref{newunif}) and
deals with properties of the matrix $(g_0)'(v_0)$ instead of the
Lipschitz constant of $g_0.$ For any set $M\subset [0,T]$
measurable in the sense of Lebesgue we denote by ${\rm mes}(M)$
the Lebesgue measure of $M.$

{\bf Theorem 2.} {\it Let $g\in
C^0(\mathbb{R}\times\Omega\times[0,1],\mathbb{R}^k)$ satisfying
(\ref{lip}). Let $g_0$ be the averaging function given by
(\ref{averfunction}) and consider $v_0\in\Omega$ such that
${g}_0(v_0)=0.$ 
Assume that
\begin{enumerate}
\setcounter{enumi}{\value{items}}
\item\label{newdef} given
 any $\widetilde{\gamma}>0$  there exist $\widetilde{\delta}>0$ and $M\subset [0,T]$
measurable in the sense of Lebesgue with ${\rm
mes}(M)<\widetilde{\gamma}$ such that for every $v\in
B_{\widetilde{\delta}}(v_0),$  $t\in[0,T]\setminus M$ and
$\eps\in[0,\widetilde{\delta}]$ we have that $g(t,\cdot,\eps)$ is
differentiable at $v$ and
$\|g'_v(t,v,\eps)-g'_v(t,v_0,0)\|\le\widetilde{\gamma}$.

 \setcounter{items}{\value{enumi}}
\end{enumerate}
Finally assume that
\begin{enumerate}
\setcounter{enumi}{\value{items}}
\item\label{NE0o}
$g_0$ is continuously differentiable in a neighborhood of $v_0$
and the real parts of all the eigenvalues of $({g}_0)'(v_0)$ are
negative. \setcounter{items}{\value{enumi}}
\end{enumerate}
Then there exists $\delta_1>0$ such that for every
$\eps\in(0,\delta_1]$,  system (\ref{ps}) has exactly one
$T$--periodic solution $x_\eps$ with $x_\eps(0)\in
B_{\delta_1}(v_0).$ Moreover the solution $x_\eps$ is
asymptotically stable and $x_\eps(0)\to v_0$ as $\eps\to 0.$}

For proving Theorem~2 we observe that the property (v) implies
(ii), while the property (vi) implies (iv). The latter observation
is totally based on the lemma~2 from the Krasnosekskii's book [5]
(P.~91).

2. Consider now an application of the result announced.

In his paper \cite{hog} Hogan first proved the existence of a
limit cycle for the nonsmooth van der Pol equation $\ddot
u+\eps(|u|-1)\dot u+u=0.$ By means of the results from the
previous section we derive now the resonance curves determining
the location of  stable periodic solutions of the perturbed
equation
\begin{equation}\label{vp}
\ddot u+\eps \left(|u|-1\right)\dot u+(1+a\eps)u=\eps\lambda\sin
t,
\end{equation}
where $a$ is a detuning parameter and $\eps\lambda\sin t$ is an
external force.

Following Andronov and Witt \cite{andr} we are concerned with the
dependence of the amplitude of $2\pi$-periodic solutions of (3) on
 $a$ and $\lambda$ providing that $\eps>0$ is sufficiently small.
 By the other words we look for such $A$ and $\phi$ which determine
 $2\pi$-periodic solutions of (3) converging as $\eps\to 0$ to
\begin{equation}\label{eq}
  u(t)=A\sin(t+\phi),\ \ \dot u(t)=A\cos(t+\phi).
\end{equation}

The assumption (v) of theorem~2 applied to (3) (after rewritting
(3) in the form (1)) turns out to be satisfied for any
$a,\lambda_0\in\mathbb{R},$ while the assumption ${\rm det}
(g_0)'((A\sin\phi,A\cos\phi))\not=0$  (required for validity of
(vi)) brings us to the following algebraic equation
\begin{equation}\label{eqA}
  A^2\left(a^2+\left(1-\frac{4}{3\pi}|A|\right)^2\right)=\lambda^2.
\end{equation}
Furthermore the negativity assumption for the real parts of
eigenvalues of $(g_0)'(M,N)$ takes the form of the following two
inequalities
\begin{equation}\label{depos}
\pi^2(1+a^2)+\frac{32}{9}(M^2+N^2)-4\pi\sqrt{M^2+N^2}>0,
\end{equation}
and
\begin{equation}\label{muneg}
  2\left(\pi-2\sqrt{M^2+N^2}\right)<0.
\end{equation}
Therefore theorem~2 allows to conclude that the solution
$u(t)=A\sin(t+\phi)$ generates asymptotically stable
$2\pi$-periodic solutions of equation (3) providing that $A,$
$M=A\sin\phi,$ $N=A\cos\phi$ satisfy (5)-(7).

We note that the corresponding assumptions for the classical Van
der Pol  equation $ \ddot u+\eps \left(u^2-1\right)\dot
  u+(1+a\eps)u=\eps\lambda\sin t$ have the following form (see the formulas~(5.21)
and (16.6)--(16.7) of Malkin's book \cite{mal})
\[
A^2\left(a^2+\left(1-\frac{A^2}{4}\right)^2\right)=\lambda^2,
\]
\[
1+a^2-(M^2+N^2)+\frac{3}{16}(M^2+N^2)^2>0,
\]
\[
2-(M^2+N^2)<0
\]
respectively.

\smallskip

The second author is partially supported by a \linebreak MEC/FEDER
grant number MTM2005-06098-C02-01 and by a CICYT grant number
2005SGR 00550 and the third author is partially supported by the
Grant BF6M10 of Russian Federation Ministry of Education and CRDF
(US), and by RFBR Grants 06-01-72552, 05-01-00100.

The authors thank Professor A. I. Perov for discussions of the
results.

\

\noindent {\bf A. Buic\u{a}:} Department of Applied Mathematics,
        Babe\c{s}--Bolyai University, Cluj--Napoca, Romania ({abuica@math.ubbcluj.ro}).
\vskip0.2cm

\noindent {\bf J.~Llibre:}  Departament de Matem\`atiques,
        Universitat Aut\`onoma de Barcelona, 08193
        Bellaterra, Barcelona, Spain ({jllibre@mat.uab.cat}).

\vskip0.2cm \noindent {\bf O. Yu. Makarenkov (corresponding
author):} Research Institute of Mathematics, Voronezh State
University, 394006, Universitetskaja pl. 1, Voronezh, Russia
({omakarenkov@math.vsu.ru}).

\end{document}